\documentclass[draft]{amsart}

\usepackage{pgf,tikz, ytableau}
\usetikzlibrary{arrows}
\usepackage{cite}

\usepackage{amsmath}
\usepackage{amssymb,graphicx, enumerate, color, hyperref}

\numberwithin{equation}{section}
 \theoremstyle{plain}
\newtheorem{theorem}{Theorem}[section]

\theoremstyle{remark}

\theoremstyle{empty}

\begin{document}

\makeatletter
\def\imod#1{\allowbreak\mkern10mu({\operator@font mod}\,\,#1)}
\makeatother

\author{Ali Kemal Uncu}
   \address{University of Bath, Faculty of Science, Department of Computer Science, Bath, BA2 7AY, UK}
   \email{aku21@bath.ac.uk}

\title{A Weighted Words Study of MacMahon's and Russell's Modulo 6 Identities}

\dedicatory{To Krishnaswami Alladi, on the occasion of his 70th birthday}

\begin{abstract} We give new proofs of MacMahon and Russell's modulo 6 identities using the method of weighted words. We also present a new refinement of MacMahon's identity, some related finite sum identities, and a companion partition theorem to sequence avoiding partitions theorem of the author and Andrews.
\end{abstract}   

  

\date{\today}

\maketitle

\section{Introduction}

 In 1993, Alladi and Gordon introduced the method of weighted words \cite{Generalized_Schur}. While studying Schur's partition theorem, by assigning weights (or colors) to different parts of partitions with the gap conditions, they generalized the problem and later found key identities that led to the proof. This approach has proven to be successful in various partition theory problems \cite[etc.]{Alladi, Alladi_Gollnitz, Alladi_Gollnitz1, Refinement, Weighted_RR, FourPar, Generalized_Schur}. However, discovering key identities is highly non-trivial, and it has been among the limiting factors of the original method. Recently, Dousse~\cite[etc.]{Dousse_Unif, Dousse_Siladic, Dousse_revisit, Dousse_Primc} revisited this method and changed its main line of finding the necessary analytic identities to instead form recurrences over the largest part of partitions with gap conditions and working on these recurrences and shift equations. This added flexibility to this method and increased its applicability. Following that lead, the author and Ablinger extended the method slightly and implemented the recurrence forming and manipulation in their Mathematica package \texttt{qFunctions} \cite{qFuncs}. Here, we will give proofs of and find some refinements of two partition theorems using the method of weighed words.

 Throughout this paper, we follow the definitions of \cite{Theory_of_Partitions}. A \textit{partition} is a finite non-increasing sequence of positive integers. Elements of the list that make a partition are called \textit{parts}. If the total of a partition is $n$, that partition is called \textit{a partition of }$n$. For example, the empty list is the (unique) partition of 0, and (1,1) is a partition of 2. 

In the second volume of his influential book \cite{MacMahon}, Major P. MacMahon proved the following modulo 6 identity.

\begin{theorem}[MacMahon] \label{thm_MacMahon} Let $n$ be a non-negative integer. The number of partitions of $n$ with no consecutive integers as parts and all parts $\geq 2$ is equal to the number of partitions of $n$, where no part occurs exactly once, is equal to the number of partitions of $n$ into parts congruent to $0,2,3$ or $4$ modulo 6.    
\end{theorem}

MacMahon's partition theorem has three parts: the one where the partitions satisfy some gap conditions between consecutive parts, the one where the partitions are counted by the frequencies (a.k.a number of occurrences) of parts, and the one where the parts satisfy some congruence conditions. Andrews generalized the frequency-congruence connection of this identity half a century after the original identity's discovery \cite{GeorgeMacMahon}. This was later generalized by Subbarao \cite{Subbarao}. Notably, most of the original proofs of these identities are due to generating function manipulations. Some bijective proofs of these identities recently started to appear in the literature  \cite{GeorgeHenriketal, FuSellers, BM, BM2}.

At the Legacy of Ramanujan 2024 conference \cite{Leg2024website}, where partition theory and $q$-series community come together to celebrate the 85th birthdays of Andrews and Berndt, Russell presented a new companion theorem to MacMahon's partition identity. 

\begin{theorem}[Russell]\label{thm_Russell} Let $n$ be a non-negative integer. The number of partitions of $n$ with no $2$s, where if adjacent parts differ by $1$ or $2$ then their sum is $\not\equiv 0$ or $\equiv 0$ modulo 3, respectively, is equal to the number of partitions of $n$ into parts congruent to $0,1,3$ or $5$ modulo 6.    
\end{theorem}

This result and some related refinement results later appeared in \cite{Russell}.

In both identities, Theorems~\ref{thm_MacMahon} and \ref{thm_Russell}, one set of partitions involves static gap conditions between consecutive parts of a given partition. This is the perfect setup for the method of weighted words. 

Following Russell's talk, the author asked if the method of weighted words approach had been tried and whether it could also provide some more insight into these problems. After confirming that this method had not been tried, the author decided to study these two identities through this lens. This led to another refinement of MacMahon's partition identity, and some other connections between  MacMahon's identity and the author's work with Andrews \cite{AndrewsUncu}. The aforementioned new refinement of MacMahon's theorem and an example are as follows:

\begin{theorem}\label{thm_Mod2_MM_refinement} Let $n$ and $m$ be non-negative integers. The number of partitions $n$ with no consecutive integers as parts and all parts $\geq 2$ with exactly $m$ odd parts is equal to the number of partitions of n into parts congruent to $0,2,3$ or $4$ modulo 6 with exactly $m$ odd parts.
\end{theorem}

Although there is no direct mandate on the odd parts of the partitions with the gap conditions (except for the missing parts size 1), Theorem~\ref{thm_Mod2_MM_refinement} shows that Theorem~\ref{thm_MacMahon} can be refined, and the number of odd parts in the sets of partitions with the gap and congruence conditions can also be equated. For example, let $n=18$ and $m=2$. 16 partitions satisfy these conditions in both sets of partitions of Theorem~\ref{thm_Mod2_MM_refinement}, and these partitions are as follows:
\[\begin{array}{lc}
\text{Gap Conditions:}&
\begin{array}{c} (15,3),\ (13,5),\ (12,3,3),\ (11,7),\ (11,5,2),\ (10,5,3),\\ (9,9),\ (9,7,2),\ (9,6,3),\ (9,5,2,2),\ (8,5,5),\ (7,7,4),\\
(7,7,2,2),\ (7,5,2,2),\ (6,6,3,3),\ (5,5,2,2,2,2). \end{array}  \\ \hline \\[-2ex]
\text{Congruence Conditions:}&
    \begin{array}{c} (15,3),\ (12,3,3),\ (10,3,3,2),\ (9,9),\ (9,6,3),\ (9,4,3,2),\ (9,3,2,2,2),\\ (8,4,3,3),\ (8,3,3,2,2),\ (6,6,3,3),(6,4,3,3,2),\ (6,3,3,2,2,2),\\ (4,4,4,3,3),\ (4,4,3,3,2,2),\ (4,3,3,2,2,2,2),\ (3,3,2,2,2,2,2,2).
    
    \end{array}  
\end{array}\]

By studying the recurrence systems weighted words approach generates for MacMahon's theorem, we also prove the summation formulas.

\begin{theorem}\label{thm_Main_Sum} For $n\geq1$ we have
\begin{align}\label{eq_sum1}
    \sum_{i,j\geq0} & (-1)^j q^{i^2 + 3ij + \frac{3j(3j+1)}{2}} {n-i-3j\brack i}_q {n-i-2j-1\brack j}_{q^3} = \sum_{j\geq0} q^{j(3j+1)}(q^2;q^3)_j{n\brack 3j+1}_q,\\
    \label{eq_sum2}\sum_{i,j\geq0} &(-1)^j q^{i^2 + 3ij + \frac{3j(3j+1)}{2}} {n-i-3j\brack i}_q {n-i-2j-1\brack j}_{q^3} \\\nonumber&- q\sum_{i,j\geq0} (-1)^j q^{i(i+1) + 3ij + \frac{9j(j+1)}{2}} {n-i-3j-1\brack i}_q {n-i-2j-2\brack j}_{q^3}= \sum_{j\geq0} q^{j(3j-1)}(q;q^3)_j{n\brack 3j}_q.
\end{align}
\end{theorem}

We also prove a companion theorem to the author's partition theorem with Andrews on sequence avoiding overpartitions \cite[Corollary 1.3]{AndrewsUncu}:

\begin{theorem}\label{thm_main_comp} Let $n$ be a non-negative integer
    The number of overpartitions (a partition where the first appearance of a part may be overlined) of $n$, where the first instance of $1$ is always overlined and $(\overline{k+1},k+1,\dots, k+1, \overline{k})$ is allowed, but no $(k+1, \overline{k})$ pattern appears if the overlined $k+1$ part is missing is equal to the number of partitions of $n$ into two colors red and green where the green parts are always 2 modulo 3.
\end{theorem}

The plan of the paper is as follows. We present necessary definitions and a brief motivation and introduction to the method of weighted words in Section~\ref{S_back}. We prove Theorems~\ref{thm_MacMahon} and \ref{thm_Mod2_MM_refinement} in Section~\ref{Sec_MM}. Section~\ref{Sec_R} has a proof of Theorem~\ref{thm_Russell} and its refinement Theorem~\ref{thm_Russell_refinement}. Section~\ref{Sec_RMM} presents a proof and a refinement of Theorem~\ref{thm_MacMahon} ( Theorem~\ref{thm_Mod3_MM_refinement}), which is equivalent to Russell's refinement of the very theorem. In Section~\ref{S_sums}, we present some recurrence relations and reductions to find closed formulas related to the generating functions of the partitions with gap conditions in Theorem~\ref{thm_MacMahon}, and also present Theorems~\ref{thm_Main_Sum} and \ref{thm_main_comp}. This also shows the connections between the results of this paper and \cite{AndrewsUncu}. The last section is reserved for possible plans to pursue following this research.

\section{Necessary Definitions and Some Background}\label{S_back}

Let $\pi$ be a partition; we denote the size (total of its parts) of $\pi$ by $|\pi|$. Let $A$ be a set of partitions. The \textit{generating function} for the partitions from $A$ (grouped with respect to their sizes) is given by the formal sum \[\sum_{\pi\in A}q^{|\pi|},\] commonly called the \textit{generating function for the number of partitions}. In this formal sum, once unfolded, the coefficient of some fixed exponent of $q$ in the power series yields the number of different partitions from $A$ that have that fixed size. More precisely, for any fixed $n$, the coefficient of $q^n$ in the formal sums is equal to the number of partitions in $A$ with size $n$. The generating functions can be refined by adding more statistics into consideration. Let $\nu_o(\pi)$ be the number of odd parts in partition $\pi$. Then, the sum $\sum_{\pi\in A} q^{|\pi|} a^{\nu_o(\pi)}$ is the generating function for the number of partitions as before, where in addition the exponent of $a$ keeps track of the number of odd parts in the partitions.

We define the $q$-Pochhammer symbol \[
(a;q)_n := \prod_{i=0}^{n-1} (1-a q^i),\ \ \text{and}\ \  (a_1,a_2,\dots,a_k;q)_n := \prod_{i=1}^{k} (a_i;q)_n, 
\]
where $k\in\mathbb{Z}_{>0}$ and $n \in \mathbb{Z}_{\geq 0} \cup \{\infty\}$. Also, we define the $q$-binomial coefficients as \[{a+b \brack b}_q = \left\{\begin{array}{cc}
	\frac{(q;q)_{a+b}}{(q;q)_a (q;q)_b}, & \text{if } a,b\geq 0,\\
	0& \text{otherwise.}
\end{array}\right.\] The classical partition theoretic interpretation of the $q$-binomial coefficient \[{a+b\brack b }_q\] is the generating function for the number of partitions where the largest part is $\leq a$ and the number of parts $\leq b$.

Let $P$ be the set of all partitions. It is easy to see that, the generating function for the number of partitions has a closed form formula \[\sum_{\pi\in P} q^{|\pi|} = \frac{1}{(q;q)_\infty}.\] Moreover, the infinite product converges absolutely for $|\pi|<1$. Similarly, let $P_{\{i_1,i_2,\dots,i_k\},m}$ be the set of partitions where the parts are from the residue classes $i_1,i_2,\dots,i_k$ modulo $m$ (where all $i_j$s are distinct, and if $i_j=0$ it is replaced with $i_j=m$). Then, \[\sum_{\pi\in P_{\{i_1,i_2,\dots,i_k\},m}} q^{|\pi|} = \frac{1}{(q^{i_1},q^{i_2},\dots,q^{i_k};q^m)_\infty}.\] For the bounded versions and generalizations of these generating functions, we refer to \cite{Theory_of_Partitions} and \cite{BU}.

It should now be clear that the partitions counted by the congruence conditions of Theorem~\ref{thm_MacMahon} have the generating function $1/(q^2,q^3,q^4,q^6;q^6)_\infty$. Hence, if we can prove that the generating function for the partitions counted by the gap conditions also has the same generating function, it would prove that the partitions from the two sets are equinumerous. The method of weighted words aims to do this using combinatorial reasoning. One starts by fixing a number of colors (weights), then associates the list of positive integers with the levels of the colors. For example, let's pick three colors $a,b,$ and $c$. We match the positive integers with levels of these colors as follows:
\[\begin{array}{cccccccc}
\begin{array}{c}1\\a_1\end{array} & \begin{array}{c}2\\b_1\end{array} & \begin{array}{c}3\\c_1\end{array} & \begin{array}{c}4\\a_2\end{array} & \begin{array}{c}5\\b_2\end{array} & \begin{array}{c}6\\c_2\end{array} & \begin{array}{c}7\\a_3\end{array} & \dots
\end{array}\]
This is akin to reading splitting positive integers into three residue classes and then attaching $a$, $b$, and $c$ tags to the residue classes $1$, $2$, and $0$, respectively, with a clear correspondence $(a_i,b_i,c_i)\mapsto (3i-2, 3i-1, 3i)$. One then translates the gap conditions of the theorem into gap conditions between different levels of the colors. For example, Theorem~\ref{thm_MacMahon} says that no consecutive parts appear together as parts in the partitions. Among others, this implies that ``if $b_k$ is a part of the partition, then the largest $a$ colored part that can be next is $a_{k-1}$." The collection of these conditions can be written into matrix form, which is called a \textit{transition matrix}. One then defines the \textit{weights} of partitions into colors, where the part $x_i \mapsto x q^i$. The aim of the method of weighted words then is to find a formula for the generating function of the partitions counted with respect to some transition matrix with as many colors intact, so that one can later, by using summation formulas or by shift equations, show that the generating function of the partitions counted with respect to the gap conditions is the same as the generating function as the partitions counted with respect to the congruence conditions. Interested readers are invited to explore \cite{qFuncs} for a detailed explanation of the technique, its implementation, and its slight generalization.

\section{A proof of MacMahon's theorem by 2-color weighted words}\label{Sec_MM}

One can express the gap conditions of MacMahon's theorem (Theorem~\ref{thm_MacMahon}) with 2 colors. Let \begin{equation}\label{2colorlabel}a_1 \leq b_1 \leq a_2 \leq b_2\leq \dots,\end{equation} where, on the back of our minds, $a$'s and $b$'s simulate odd and even positive integers, respectively. More precisely, the elements in the ordering \eqref{2colorlabel} contribute the weights \[aq,bq,aq^2,bq^2,\dots.\] Under the mapping $(a,b,q)\mapsto(1/q,1,q^2)$, these weights turn into the usual contributions of the parts of partitions to the generating functions.

The gap conditions of the first set of partitions described in MacMahon's theorem (Theorem~\ref{thm_MacMahon}) say that if some part appears in a partition, then more copies of this part can appear as parts. However, if a number $x$ is a part of a partition, then $x\pm 1$ are not parts of that partition. In the 2 colored alphabet of \eqref{2colorlabel}, we can encode these rules using a transition matrix. 

The transition matrix for MacMahon's theorem is as follows:
\[
M = \bordermatrix{~ & a & b \cr
              a & 0 & 2 \cr
              b & 1 & 0 \cr}.
\]
Here, the entries of the matrix $M$ should be read as the descent rules of indices from the color on the row to the color on the column. For example, the $M_{1,1}$ entry 0 is interpreted as the statement ``if $a_k$ appears as a part in a partition, then $a_{k-0}$ can also be a part." In other words, the part $a_k$ can be repeated in a partition. Similarly, the $M_{1,2}$ entry 2 should be interpreted as the statement ``if $a_k$ is a part of a partition then $b_{k-2}$ can follow it (not $b_{k-1}$)."

Let $g^M_{a_{n}}(a,b;q)$ and $g^M_{b_{n}}(a,b;q)$ be the generating functions for all partitions in 2 colors that satisfy the transition conditions of $M$, where the largest part is bounded by $a_n$ and $b_n$, respectively. We have the initial values \[g^M_{a_{1}}(a,b;q)=1 \quad\text{  and  }\quad  g^M_{b_{1}}(a,b;q) = \frac{1}{1-b q}.\] In MacMahon partitions, 1 does not appear as a part. This is reflected in the initial values as $a_1$ not being a part of the 2 colored partitions. 

Also note that \[g^M(a,b;q):= \lim_{n\rightarrow\infty} g^M_{a_{n}}(a,b;q) = \lim_{n\rightarrow\infty} g^M_{b_{n}}(a,b;q). \]

We can automatically generate the coupled recurrences satisfied by $g^M_{a_n}(a,b;q)$ and $g^M_{b_n}(a,b;q)$ using \texttt{qFunctions} \cite{qFuncs}. These recurrences are
\begin{align}
\label{eq:ga_rec}g^M_{a_n}(a,b;q) =&   \frac{a q^n }{1 - a q^n} g^M_{a_{n-1}}(a,b;q) + g^M_{b_{n-1}}(a,b;q),\\
\label{eq:gb_rec}g^M_{b_n}(a,b;q) =& g^M_{a_n}(a,b;q) + \frac{b q^n}{1 - b q^n} g^M_{b_{n-1}}(a,b;q).
\end{align}

These recurrences can be used to generate more initial terms of the sequences and be automatically uncoupled by \texttt{qFunctions} \cite{qFuncs}, using \cite{HolonomicFunctions} under the hood. These uncoupled recurrences are as follows:

\begin{align} 
   \nonumber (1 - a q^{n+1}) (1 &- b q^{n+1}) (1 - a q^{n+2}) g^M_{a_{n+2}}(a,b;q) \\\label{unocup_rec1_MM}&- (1 - a q^{ n+1}) (1 - a b q^{ 2 n+3}) g^M_{a_{n+1}}(a,b;q) +a b q^{2n + 2} (1 - a q^{ n+2}) g^M_{a_{n}}(a,b;q)  =0\\
      \nonumber (1 - b q^{n+1}) (1 &- a q^{n+2}) (1 - b q^{n+2}) g^M_{b_{n+2}}(a,b;q) \\\label{unocup_rec2_MM}&- (1 - b q^{ n+1}) (1 - a b q^{ 2 n+4}) g^M_{b_{n+1}}(a,b;q) +a b q^{2n + 3} (1 - b q^{ n+2}) g^M_{b_{n}}(a,b;q)  =0
\end{align}

These recurrences are visually similar. This is common in the method of weighted words. To prove MacMahon's theorem, we now look for a shift equation. This is Dousse's spin on the method of weighted words \cite{Dousse_revisit}. The more free parameters we can keep while identifying the shift equation, the better. By looking at the initial conditions of the two generating function sequences  $g^M_{a_{n}}(a,b;q)$ and $g^M_{b_{n}}(a,b;q)$, we can observe that
\begin{equation}\label{MacMahon_shift}
    g^M_{a_{n+2}}(a,1;q) = \frac{1-a q^{n+3}}{(1-q^{n+1})(1-a q^2)} g^M_{b_n}(a q^3,1;q)
\end{equation}
gets satisfied. We can prove this observation by making the needed substitutions in the uncoupled recurrences \eqref{unocup_rec1_MM} and \eqref{unocup_rec2_MM}. The index shift $n\mapsto n+2$ and the variable substitution $b\mapsto 1$ in \eqref{unocup_rec1_MM}, and the formal substitution of the leading coefficient following the substitution $(a,b)\mapsto(aq^3,1)$ in \eqref{unocup_rec2_MM} shows that these two sequences satisfy the recurrence
\begin{align*}(1 - q^{n+3}) (1 &- aq^{n+3}) (1 - a q^{n+4}) c_{n+2} - (1 - a q^{ n+3}) (1 - a q^{ 2 n+7}) c_{n+1} +a  q^{2n + 6} (1 - a q^{ n+4}) c_n  =0. \end{align*} Hence, showing that the first two initial values of the left- and right-hand sides are equal is enough to prove \eqref{MacMahon_shift}.

The equation \eqref{MacMahon_shift} with the convergence condition $|q|<1$ takes us to the doorstep of the product generating function of Theorem~\ref{thm_MacMahon} (which was presented in Section~\ref{S_back}). First, we take the limit $n\rightarrow \infty$. With the convergence condition, this yields
\begin{align} \label{MM_pure_shift} g^M(a,1;q) &= \frac{1}{1-a q^2} g^M(a q^3,1;q), \\
\intertext{which can be used iteratively on itself:}
\nonumber&= \frac{1}{(1-a q^2)(1-aq^5)} g^M(a q^6,1;q), \\
\nonumber&= \frac{1}{(1-a q^2)(1-aq^5)(1-aq^8)} g^M(a q^9,1;q), \\
\nonumber&\vdots\\
\label{MM_proof_last_step}&= \frac{1}{(a q^2 ;q^3)_\infty}g^M(0,1;q).
\end{align}

Each application of \eqref{MM_pure_shift} on the right-hand side of the equation introduces one more $q$-factor to the coefficient, and also multiplies the first argument of the generating function with $q^3$. In the limit $n\rightarrow\infty$, this turns into the $q$-Pochhammer term in \eqref{MM_proof_last_step} and the first argument of the generating function vanishes by imposing the convergence condition $|q|<1$. 

The final step is to observe that $g^M(0,1,q)$ is the generating function for the number of all partitions. Plugging in $a=0$ makes all the partitions that satisfy the transition rules of $M$ irrelevant. Only the partitions into $b$ colors are counted by $g^M(0,1,q)$. Moreover, the transition rules say that $b_n$ terms can repeat and a $b_n$ part can be followed by a $b_{n-1}$ part. This is all the partitions. Hence, we get that \[ g^M(0,b;q) = \frac{1}{(b q ; q)_\infty}.\]

The above paragraph and \eqref{MM_proof_last_step} imply the following theorem.

\begin{theorem}\label{Weighted_MM} For $|q|<1$, we have 
    \[g^M(a,1;q) = \frac{1}{(a q^2 ;q^3)_\infty(q;q)_\infty},\]
where $g^M$ is the generating function for the number of partitions counted with the gap conditions encoded by $M$, the exponent of $q$ counts the total size of the partitions, and the exponent of $a$ keeps track of the number of $a$ colored parts in the partitions.
\end{theorem}

As mentioned at the beginning of the section, a direct consequence of this theorem is MacMahon's original theorem. Letting $(a,q)\mapsto (1/q, q^2)$ yields Theorem~\ref{thm_MacMahon}. The left-hand side generating function of the partitions that satisfy the gap conditions prescribed in Theorem~\ref{thm_MacMahon} and the right-hand side product after substitutions turns into \begin{equation}\label{MM_prod_GF}\frac{1}{(q^2,q^3,q^4,q^6;q^6)_\infty},\end{equation} which is the generating function for the number of partitions into parts congruent to $0,2,3$ or $4$ modulo 6.

Moreover, Theorem~\ref{Weighted_MM}, with the presence of the extra parameter $a$, allows us to refine this result by keeping a count on odd parts after the necessary substitutions, which implies Theorem~\ref{thm_Mod2_MM_refinement}. The proof of Theorem~\ref{thm_Mod2_MM_refinement} comes from the substitution $(a,q)\mapsto(a/q,q^2)$ in Theorem~\ref{thm_MacMahon} followed by comparing the coefficients of the exponents of $a$ in the formal power series.

\section{A proof of Russell's Theorem by 3-color weighted words}\label{Sec_R}

One needs 3 colors to express the gap conditions of Russell's Theorem (Theorem~\ref{thm_Russell}). Let \begin{equation}\label{3color}a_1 \leq b_1 \leq c_1 \leq a_2 \leq b_2 \leq c_2 \leq \dots\end{equation} be the three colors and their ordering. We take the following transition matrix 
\[R:= \bordermatrix{~ & a & b & c\cr
              a & 0 & 1 & 1\cr
              b & 1 & 0 & 2\cr
              c & 1 & 0 & 0},
\]which turns to the gap conditions of Theorem~\ref{thm_Russell} under the mapping $(a,b,c,q)\mapsto(1/q^2, 1/q,1,q^3)$ that takes the $a, b$ and $c$ colored parts to $1,2$ and $0$ modulo 3 parts, respectively.

Let $g^R_{a_{n}}(a,b,c;q)$, $g^R_{b_{n}}(a,b,c;q)$ and $g^R_{c_{n}}(a,b,c;q)$ be the generating functions for all partitions in 3 colors that satisfy the transition conditions of $R$, where the largest part is $a_n$, $b_n$ and $c_n$, respectively. We have the initial values \begin{equation}\label{R_initial}g^R_{a_{1}}(a,b,c;q)= g^R_{b_{1}}(a,b,c;q)=\frac{1}{1-a q} \quad\text{  and  }\quad  g^R_{c_{1}}(a,b,c;q) =1+ \frac{ a q}{1-a q} +\frac{ c q}{1-c q}.\end{equation} The term $b_1$ is missing in the partitions, which corresponds to ``no 2's as parts" condition of Theorem~\ref{thm_Russell}.

These generating functions satisfy the following coupled system of recurrence relations
\begin{align*}
    g^R_{a_{n}}(a,b,c;q) &= \frac{1}{1- a q^n} g^R_{c_{n-1}}(a,b,c;q),\\
    g^R_{b_{n}}(a,b,c;q) &= g^R_{a_{n}}(a,b,c;q) + \frac{b q^n}{1- b q^n} g^R_{b_{n-1}}(a,b,c;q),\\
    g^R_{c_{n}}(a,b,c;q) &= g^R_{b_{n}}(a,b,c;q) + \frac{b c q^{2n}}{(1- b q^n)(1- c q^n)} g^R_{b_{n-1}}(a,b,c;q) + \frac{c q^n}{1- c q^n} g^R_{c_{n-1}}(a,b,c;q).
\end{align*}

Similar to Section~\ref{Sec_MM}, we can uncouple the recurrences automatically. One of which is the following:
\begin{align}
\nonumber
(1 -b q^{n+1}) &(1- c q^{n+1}) (1 -a q^{n+2}) g^R_{a_{n+2}}(a,b,c;q)\\
    \label{R_an_rec}&- (1 - b c q^{ 2 n+1} - a b q^{2 n+2} - a c q^{ 2 n+2} + a b c q^{3 n+2} + a b c q^{3 n+3}) g^R_{a_{n+1}}(a,b,c;q) \\
    \nonumber&+a b c q^{3 n+1}g^R_{a_{n}}(a,b,c;q) =0
\end{align}

We now experiment and observe a shift equation between these generating functions while keeping as many free parameters in the problem. This search yielded the relation
\begin{equation}\label{Russell_shift}
    g^R_{a_{n+1}} (a,b,1;q) = \frac{1- b q^{n+2}}{(1-b q^2)(1-a q^{n+1})} g^R_{b_{n}}(a,b q^2,1;q)
\end{equation}
Under the assumption $|q|<1$, after $n\rightarrow\infty$, we get a shift equation that can be used iteratively from \eqref{Russell_shift}:
\begin{equation}\label{Russell_Pure_Shift}
    g^R(a,b,1;q) = \frac{1}{(1-b q^2)} g^R(a,b q^2,1;q),
\end{equation}
where $g^R(a,b,c;q) := \lim_{n\rightarrow\infty}g^R_{a_{n}}(a,b,c;q)= \lim_{n\rightarrow\infty}g^R_{b_{n}}(a,b,c;q)= \lim_{n\rightarrow\infty}g^R_{c_{n}}(a,b,c;q).$ Using \eqref{Russell_Pure_Shift} repeatedly on the right-hand side generating function gives \begin{equation}\label{Russell_Pure_Shift2}
    g^R(a,b,1;q) = \frac{1}{(b q^2;q^2)_\infty} g^R(a,0,1;q).
\end{equation}
Now, we focus on finding a closed formula for the function $g^R(a,0,1;q)$. Substituting $(a,b,c,q)\mapsto (a,0,1,q)$ in \eqref{R_an_rec}, we get the order one recurrence \[ g^R_{a_{n+1}} (a,0,1,q) = \frac{(1 - a q^{2n})}{(1 - q^{n} )(1 - a q^{n+1})} g^R_{a_{n}} (a,0,1,q) \]
Unrolling this recurrence relation gives \[ g^R_{a_{n+1}} (a,0,1,q) = \frac{(a q^{2};q^2)_n}{(q;q)_n (a q^2;q)_n} g^R_{a_{1}} (a,0,1,q).\] Using the initial condition \eqref{R_initial} with the necessary substitutions (although they do not change the outcome in this example), we get
\[g^R_{a_{n+1}} (a,0,1,q) =\frac{(a q^{2};q^2)_n}{(q;q)_n (a q;q)_{n+1}}. \]
Hence, \[g^R(a,0,1,q) =\frac{1}{(q;q)_\infty (a q;q^2)_\infty}.\] We get the following theorem by plugging in $g^R_{a_{n+1}} (a,0,1,q)$ in \eqref{Russell_Pure_Shift2}.

\begin{theorem}\label{Weighted_R} For $|q|<1$, we have 
    \[g^R(a,b,1;q) = \frac{1}{(q;q)_\infty (a q;q^2)_\infty(b q^2;q^2)_\infty},\]
where $g^R$ is the generating function for the number of partitions counted with the gap conditions encoded by $R$, the exponent of $q$ counts the total size of the partitions, and the exponents of $a$ and $b$ keep track of the number of $a$ and $b$ colored parts in the partitions, respectively.
\end{theorem}

After the substitutions $(a,b,q)\mapsto(1/q^2,1/q,q^3)$, Theorem~\ref{Weighted_R} proves Theorem~\ref{thm_Russell}. Furthermore, the more general substitution, $(a,b,q)\mapsto(a/q^2,b/q,q^3)$ yields Russell's refinement of Theorem~\ref{thm_Russell} by comparing the exponents of $a$ and $b$ in the generating functions \cite[Theorem 9]{Russell}:

\begin{theorem}\label{thm_Russell_refinement} Let $n, m$ and $k$ be non-negative integers. The number of partitions $n$ that are counted by the gap conditions in Theorem~\ref{thm_Russell} with the added conditions that it has exactly $m$ and $k$  1 and 2 modulo 3 parts, respectively, is equal to the number of partitions of n into parts congruent to $0,1,3$ or $5$ modulo 6 with exactly $m$ 1 modulo 6 parts and $k$ 5 modulo 6 parts.
\end{theorem}

It should be noted that Theorem~\ref{Weighted_R} and \cite[Theorem 8]{Russell} that imply Theorem~\ref{thm_Russell_refinement} differ slightly. This is clear from the product sides, where Russell's $q$-Pochhammer product has $(bq;q^2)_\infty$ in the place of $(b q^2;q^2)_\infty$ of Theorem~\ref{Weighted_R}. This is due to the unnatural number ordering Russell worked with. This order is also in the scope of weighted words. Alladi \cite{Alladi} studied the impact of such different orders in the scope of weighted works in his earlier work. The transition matrix and the substitutions that would yield \cite[Theorem 8]{Russell} are \[R':= \bordermatrix{~ & a & b & c\cr
              a & 0 & 1 & 2\cr
              b & 1 & 0 & 1\cr
              c & 0 & 1 & 0}\] 
and $(a,b,c,q)\mapsto(1/q^2, 1, q^2 ,q^3)$. This substitution takes the $a,b$ and $c$ colored numbers to $1,0$ and $2$ modulo 3 parts, respectively, where the $2$ modulo 3 parts start from 5.

\section{Russell's Refinement of MacMahon's Theorem}\label{Sec_RMM}

Russell found a refinement of Theorem~\ref{thm_MacMahon} similar to Theorem~\ref{thm_Russell_refinement}. We can also give a weighted word proof of that using the transition matrix $R'$. 

Observe that the rules $R'$ impose on the 3 color order \eqref{3color} collapse to nothing but ``no consecutive parts". Hence, by the mapping $(a,b,c,q)\mapsto(1/q^2,1/q,1,q^3)$ and the necessary initial condition that all parts~$\geq 2$, we will recover the gap conditions of MacMahon's theorem (Theorem~\ref{thm_MacMahon}).

To that end, we let $g^{R'}_{a_{n}}(a,b,c;q)$, $g^{R'}_{b_{n}}(a,b,c;q)$ and $g^{R'}_{c_{n}}(a,b,c;q)$ be the generating functions for all partitions in 3 colors that satisfy the transition conditions of $R'$, where the largest part is $a_n$, $b_n$ and $c_n$, respectively. We have the initial values \begin{equation}\label{R2_initial}g^{R'}_{a_{1}}(a,b,c;q)=1, \quad  g^{R'}_{b_{1}}(a,b,c;q)=\frac{1}{1-b q} \quad\text{  and  }\quad  g^{R'}_{c_{1}}(a,b,c;q) =1+ \frac{ b q}{1-b q} +\frac{ c q}{1-c q}.\end{equation}

Once the recurrence relations and necessary initial conditions are studied, we see that \[g^{R'}_{c_n}(a,b,1,q) = \frac{1-b q^{n+1}}{(1-q^n)(1- b q)} g^{R'}_{a_n}(a,b q^2,1,q)\] is satisfied. Following that, for $|q|<1$, we can prove that \[ g^{R'}(a,0,1,q) = \frac{1}{(q;q)_\infty (a q^2;q^2)_\infty}. \] Combining these results following the steps in Section~\ref{Sec_R} yields the following theorem.

\begin{theorem}\label{Weighted_R2} For $|q|<1$, we have 
    \[g^{R'}(a,b,1;q) = \frac{1}{(q;q)_\infty (a q^2;q^2)_\infty(b q;q^2)_\infty},\]
where $g^{R'}$ is the generating function for the number of partitions counted with the gap conditions encoded by $R'$, the exponent of $q$ counts the total size of the partitions, and the exponents of $a$ and $b$ keep track of the number of $a$ and $b$ colored parts in the partitions, respectively.
\end{theorem}

When we substitute $(a,b,c,q)\mapsto(1/q^2,1/q,1,q^3)$, this theorem implies Theorem~\ref{thm_MacMahon}. Under the mapping $(a,b,c,q)\mapsto(C/q^2,A/q,1,q^3)$, Theorem~\ref{Weighted_R2} turns into Russell's refinement of the MacMahon's Theorem \cite[Theorem 7]{Russell}:

\begin{theorem}\label{thm_Mod3_MM_refinement} Let $n, m$ and $k$ be non-negative integers. The number of partitions $n$ that are counted by the gap conditions in Theorem~\ref{thm_MacMahon} with the added conditions that it has exactly $m$ and $k$  1 and 2 modulo 3 parts, respectively, is equal to the number of partitions of n into parts congruent to $0,2,3$ or $4$ modulo 6 with exactly $m$ 4 modulo 6 parts and $k$ 2 modulo 6 parts.
\end{theorem}

The theorem follows simply by comparing the exponents of $C$ and $A$ after the mentioned substitutions.

Notice that the refinement given by Theorem~\ref{thm_Mod2_MM_refinement} is different than Theorem~\ref{thm_Mod3_MM_refinement}. On the congruence sides, one refinement fixes the number of 3 modulo 6 parts, while the other one fixes the number of 2 and 4 modulo 6 parts. 

\section{Some Sum Representations and Connections with Sequences in Overpartitions}\label{S_sums}

There are no known analytic sum formulas that can be interpreted as the generating function for the number of partitions with the gap conditions from the Theorems~\ref{thm_MacMahon} and \ref{thm_Russell} or their refinements. It is desirable to find such formulas to see the analytic counterpart and implications of these theorems. Although we cannot present these generating functions here, we will show some direct connections of MacMahon's theorem and a theorem of Andrews and the author's on sequence avoiding overpartitions~\cite{AndrewsUncu}.

Identifying a $q$-hypergeometric sum formula for $g^M_{a_n}(a,b;q)$ (see Section~\ref{Sec_MM}) requires us to find an expression that satisfy \eqref{unocup_rec1_MM}. This is known to be an infamously hard problem, sometimes referred to as the \textit{inverse Zeilberger problem}. The author himself \cite{Uncu2, Uncu3}, with Ablinger \cite{qFuncs}, and with Jim\'enez-Pastor \cite{FactorialBasis} made limited attempts to this problem through combinatorics and symbolic computation.

Nevertheless, looking at the recurrences and doing formal manipulations (which can be automatically carried out by \cite{qFuncs}) can be helpful in this task. For example, \eqref{unocup_rec1_MM} suggests that a natural substitution to simplify the recurrence is \begin{equation}\label{subs}g^M_{a_n}(a,b;q) = \frac{h^M_{n}(a,b;q)}{(a q^2; q)_{n-1}}.\end{equation} The sequence $h^M_{n}(a,b;q)$ satisfies the recurrence \begin{equation}\label{Raw_subs_rec}(1-bq^{n-1})h^M_{n}(a,b;q) = (1-a b q^{2n-1})h^M_{n-1}(a,b;q)-a b q^{2n-2}(1-a q^n)h^M_{n-2}(a,b;q),\end{equation} with the initial conditions $h^M_0(a,b;q)=0$ and $h^M_1(a,b;q)=1$. If a closed-form sum expression can be found for $h^M_n(a,b;q)$, we can relate it to Theorem~\ref{thm_MacMahon} and the product generating function \eqref{MM_prod_GF}.

Notice that \eqref{Raw_subs_rec} simplifies if one takes $b\mapsto a q$. Let $h_n(a;q) := h^M_n(a,aq;q)$, then $h_n(a;q)$ satisfies \begin{equation}\label{Subs_rec}h_{n}(a;q) = (1+a q^n)h^M_{n-1}(a;q)-a^2 q^{2n-1}h_{n-2}(a;q).\end{equation} We still have $h_0(a;q)=0$ and $h_1(a;q)=1$. This substitution takes us away from the full generality of MacMahon's theorem, but at least we can find closed-form solutions to \eqref{Subs_rec} that satisfy the initial conditions in two ways:

\begin{theorem}\label{thm_Determ} For $n\geq 2$, $h_n(a;q)$ is given as the determinant of the $(n-1)\times (n-1)$ matrix:
    \[h_n(a;q) = \det\left(\begin{array}{cccccc}
    1+a q^2 & a^2 q^{5} & 0 & 0& \dots & 0\\
    1 & 1+ a q^3 & a^2 q^{7} & 0 & \ddots &0\\
    0& 1 & 1+ a q^4 & a^2 q^{9} & \ddots & \vdots\\
    0& 0 & 1 & \ddots & \ddots & 0\\
    0& \ddots & \ddots & \ddots & \ddots & a^2 q^{2n-1}\\
    0& \ddots & \ddots & \ddots & 1 & 1+a q^{n}\\
    \end{array}\right).\]
\end{theorem}

The proof quickly follows by evaluating the determinant by expanding it with respect to its last column. 

\begin{theorem}\label{thm_Fermionic} Let $n$ be an integer. We have
    \begin{equation}
        \label{FermionicA} h_n(a;q) = \sum_{i,j\geq0} (-1)^j a^{i+3j} q^{i(i+1) + 3ij + \frac{9j(j+1)}{2}} {n-i-3j\brack i}_q {n-i-2j-1\brack j}_{q^3}.
    \end{equation}
\end{theorem}

One can easily show that this sum satisfies \eqref{Subs_rec} by automated means, such as \cite{HolonomicFunctions}.

Let $h'_n(a;q):=(a q;q)_{n}g^M_{b_n}(a,a;q)$. Then, $h'_n(a;q)$ also satisfies the recurrence \eqref{Subs_rec} (with $h_n$ replaced by $h'_n$) and it has initial values $h_0(a;q)=h_1(a;q)=1$. This implies the following theorem.

\begin{theorem} Let $n$ be an integer
    \begin{align}
        \label{FermionicB}  h'_n(a;q) &= h_n(a;q) - a^2 q^3 h_{n-1}(aq;q).
    \end{align}
\end{theorem}

Moreover, for some specific substitutions of $a$, it is possible to discover alternative single sum representations of these functions using \cite{qFuncs}. One such example is Theorem~\ref{thm_Main_Sum}. The left-hand sides of Theorem~\ref{thm_Main_Sum} are $h_n(1/q;q)$ and $h'_n(1/q;q)$, and the proofs of the identities are standard using \cite{HolonomicFunctions}. Similarly, for $h'_n(1;q)$, we have the following theorem.
\begin{theorem}\label{thm_sum_big} For $n\geq 1$, we have 
\begin{align}
    \nonumber\sum_{i,j\geq0} &(-1)^j q^{i(i+1) + 3ij + \frac{9j(j+1)}{2}} {n-i-3j\brack i}_q {n-i-2j-1\brack j}_{q^3} \\\label{eq_sum3}&-\sum_{i,j\geq0} (-1)^j q^{i(i+2) + 3ij + \frac{9j(j+1)}{2}+3j+3} {n-i-3j-1\brack i}_q {n-i-2j-2\brack j}_{q^3}\\\nonumber &= \sum_{j\geq0} q^{j(3j-1)}(q;q^3)_j{n\brack 3j-1}_q+\sum_{j\geq0} q^{j(3j+2)}(q;q^3)_j{n\brack 3j}_q.
\end{align}
\end{theorem}

We now want to take a step back and highlight the combinatorial significance of the $b\mapsto a q$ map. Section~\ref{Sec_MM} weighted words transition matrix and the initial conditions ensured that no $a$-colored parts can be followed by a $b$-colored part of one level below. Moreover, the map $(a,b,q)\mapsto(1/q,1,q^2)$ was taking these colored parts to ordinary integers. By taking $b\mapsto a q$ first, we mapped each $b_i$ to a new copy of $a_{i+1}$. In a partition, after this weighting, an $a_i$ can only come from one source: the original $a_i$ or $b_{i-1}$'s image. We can mark the first appearance of such an $a_i$ by overlining it, and any other copies of $a_i$ in the presence of $\overline{a_1}$ would be understood to be coming from the same source (as the gap conditions do not allow $a_i$ and $b_{i-1}$ to come together in the first place). Recall that $a_1$ does not appear as a part in MacMahon partitions originally. Now, the only parts that can appear in the colored setup of the MacMahon partitions are $a'_2, a_2, a'_3, a_3, \dots$, where the $a_i'$ colors are images of $b_{i-1}$'s under $b\mapsto a q$ image. Given a partition in this coloring, we can start by overlining the parts from some point in the partitions. Starting to mark parts at some part size and by alternating in colors, we overline each level going down. For example, if we start overlining parts by overlining an $a'_i$, then the next smaller part size to be overlined is $a_{i-1}$, then $a'_{i-2}$, then $a_{i-3}$, etc. Under the rules on these parts are that $\overline{a_{k+1}}$ cannot be together with $\overline{a_{k}}$ (since this can only happen either when we originally have $b_{k+1}$ followed by $a_{k+1}$ or when we have $a_{k+1}$ followed by a $b_k$), and $(\overline{a_{i+2}}, a_{i+1},\overline{a_{i}})$ cannot appear (since this configuration cannot appear in our overlining scheme). For $a\mapsto 1/q$, all parts in the partitions start at size 1, and the conditions mentioned coincide with those of \cite[Corollary 1.3]{AndrewsUncu} for overpartitions. Hence, we learn that the author and Andrews' result on sequence avoiding partition theorem is closely related to MacMahon's theorem through the 2 colored weighted words study. 

Furthermore, we can recover \cite[Theorem 1.1]{AndrewsUncu} and \cite[Corollary 1.3]{AndrewsUncu} completely by the identities we discovered already. Equation~\eqref{eq_sum1} is a new finite version of \cite[Theorem 1.1]{AndrewsUncu} and implies it in the limit $n\rightarrow \infty$. Then by Theorem~\ref{thm_Fermionic}, \eqref{eq_sum1}, \eqref{subs} with $n\rightarrow\infty$, and the interpretation made in the last paragraph, we see that \[ g^M(1/q,1;q) = \frac{1}{(q;q)_\infty (q;q^3)_\infty}.\] Here, the left-hand side is the generating function for the number of overpartitions where the pairs $(\overline{k+1},\bar{k})$ and $(\overline{k+2},k+1,\overline{k})$ do not appear as parts together, and the right-hand side is the generating function for the number of partitions into two colors red and green where each green part is congruent to 1 modulo 3. One can also see that \cite[Corollary 2.2]{AndrewsUncu} and \eqref{Subs_rec} are related up to a simple substitution and re-indexing. The same observation is true for Theorem~\ref{thm_Determ} and \cite[Corollary 2.5]{AndrewsUncu}.

Finally, in the overpartitions interpretation, plugging in $a=1$ in Theorem~\ref{Weighted_MM} gives the companion theorem to \cite[Corollary 1.3]{AndrewsUncu}. In this setup, we overline the first instance of each $b$ colored part. The part $a_1$ does not appear as a part in MacMahon's theorem, so the only part size 1 to appear in the final partitions will come from $b_1$. Hence, the first appearance of a 1 will always be overlined in our partitions. Then, the transition rules mandate that $(k+1,\dots, k+1, \overline{k})$ cannot appear (but $(\overline{k+1},k+1,\dots, k+1, \overline{k})$ is allowed). This yields Theorem~\ref{thm_main_comp}.

Here are the 12 partitions related to Theorem~\ref{thm_main_comp} when $n=5$.

\[\begin{array}{lc}
\text{Pattern Conditions:}&
\begin{array}{c} (\overline{5}),\ (5),\ (\overline{4},\overline{1}),\ (4,\overline{1}),\ (\overline{3},\overline{1},1),\ (3,\overline{1},1),\ (\overline{3},\overline{2}),\ (\overline{3},{2}),\ ({3},{2}),\\ (\overline{2},2,\overline{1}),\ (\overline{2},\overline{1},1,1),\ (\overline{1},1,1,1,1). \end{array}  \\ \hline \\[-2ex]
\text{Colored Conditions:}&
    \begin{array}{c} (5_r),\ (5_g),\ (4_r,1_r),\ (3_r,2_r),\ (3_r,2_g),(3_r,1_r,1_r),\ (2_r,2_r,1_r),\ (2_g,2_r,1_r),\\ (2_g,2_g,1_r),\ (2_r,1_r,1_r,1_r), \ (2_g,1_r,1_r,1_r),\ (1_r,1_r,1_r,1_r,1_r).
    \end{array}    
\end{array}\]

Notice that Theorem~\ref{thm_sum_big} is the finite analog of the analytic version of Theorem~\ref{thm_main_comp}. After getting rid of frivolous 0 parts, we can also interpret \eqref{eq_sum2} as a finite analog of Theorem~\ref{thm_main_comp}.

\section{Outlook}

Studying MacMahon and Russell's modulo 6 identities using weighted words not only allowed us to find a new refinement of MacMahon's identity, but also allowed us to connect this identity with some recent results of the author with Andrews. We found some new polynomial analogs of the double sum equals infinite product type identities that were discovered in \cite{AndrewsUncu}. We plan to study Russell's identity as well as the 3-color proof of MacMahon's identity in more detail in hopes of giving an independent proof of \cite[Conjecture 1.2]{AndrewsUncu}. 

There is also the question of finding a unifying companion for MacMahon's theorem. Theorems~\ref{thm_Mod2_MM_refinement} and \ref{thm_Mod3_MM_refinement} come from different-sized transition matrices. The number of colors in both approaches divides the modulus of the congruence conditions of MacMahon's theorem. We plan to attempt this theorem once again with a 6-color system to see if we can prove it and if a unifying refinement can be found through that.

The sum formulas that can be interpreted as the gap conditions of MacMahon and Russell's theorems are still missing. We aim to study these theorems more. One possible approach now is to revert to the original method of weighted words and to look for a key identity. Another possible way is to find a combinatorial construction, such as the one in \cite{Uncu2}.

It was also recently brought up to our attention that Kanade--Nandi--Russell embedded the MacMahon's identity in multiple infinite families of identities. It would be intriguing to check these partition identities and see if it would be possible to prove and refine these using the techniques mentioned here.

\section{Acknowledgement}

The author would like to thank Thieu Vo for his interest in the project and Runqiao Li for his comments.

The author would like to thank the Austrian Science Fund FWF. The research P34501N "Partition Identities Through the Method of Weighted Words" has been the main source of the many papers the author contributed to in the near past. This manuscript is among the last projects, where in its beginning, the author was still being officially partially supported by P34501N. We thank the Austrian Science Fund FWF for the great support and the trust they put in the author's abilities throughout these years.

\end{document}